\documentclass[12pt]{amsart}

\usepackage{common}
\usepackage{mathrsfs}

\title{Morley sequences in dependent theories}

\author{Alexander Usvyatsov}

\address{ Alexander Usvyatsov\\ Universidade de Lisboa \\
  Centro de Matem\'{a}tica e Aplica\cc{c}\~{o}es Fundamentais\\
  Av. Prof. Gama Pinto,2\\
  1649-003 Lisboa \\
  Portugal}

\address{Alexander Usvyatsov\\
  University of California -- Los Angeles\\
  Mathematics Department\\
  Box 951555\\
  Los Angeles, CA 90095-1555\\
  USA}

\urladdr{http://www.math.ucla.edu/\textasciitilde alexus}

\thanks{Research partially supported by Funda\cc{c}\~{a}o para a Ci\^{e}ncia e a Tecnologia, Financiamento Base 2008
- ISFL/1/209. }

\date{\today}

\begin{document}

\begin{abstract}
    We characterize nonforking (Morley) sequences in dependent theories in terms of a generalization of Poizat's
    special sequences and show that average types of Morley sequences are stationary
    over their domains. We characterize generically stable types
    in terms of the structure of the ``eventual'' type. We then
    study basic properties of ``strict Morley sequences'', based on
    Shelah's notion of strict nonforking. In particular we prove
    ``Kim's lemma'' for such sequences, and a weak version of local
    character.
\end{abstract}

\maketitle

\section{Introduction and Preliminaries}
\subsection{Introduction}
    This paper is a natural continuation of \cite{Us}, where we
    proved several useful facts on
    nonforking sequences in dependent theories which are also
    indiscernible sets. In particular, it was shown that given a
    nonforking indiscernible set $I$ over a set $A$, the global
    average type $\Av(I,\fC)$ does not fork over $A$ (and furthermore,
    if $A = \acl(A)$, then $\Av(I,\fC)$ is the unique nonforking
    extension of its restriction to $A$). These properties of
    nonforking indiscernible were crucial for understanding
    generically stable types. A natural question was: what if one
    does not assume that $I$ is an indiscernible set? Are the
    results in \cite{Us} a particular case of a general theory of
    (arbitrary)nonforking sequences in a dependent theory?

    This paper provides a complete and satisfactory answer to the
    question above. We show that a nonforking sequence (which we
    call a ``Morley sequence'' here) is in many ways the ``correct''
    object to work with in a dependent theory. Let us explain what
    we mean.

    In a stable theory a type over an algebraically closed
    set determines a unique global nonforking extension. This is not
    the case for an arbitrary type in a dependent theory, the
    simplest example being the type $p$ ``at infinity'' over the model
    $\setQ$ in the theory of dense linear orders without endpoints.
    This type has two global nonforking extensions: one is the type
    of the cut $\setQ^+$ which is finitely satisfiable in $\setQ$,
    and another one is the type ``at infinity'' over the monster
    model (which is the extension with respect to the natural
    definition of $p$). One of the main results of this paper is
    that whereas a type $p$ does not determine a unique global nonforking
    extension, such an extension \emph{is} determined by a Morley
    sequence in $p$. More precisely, we show that given a Morley
    sequence $I$ over a set $A$, the type $\Av(I,A\cup I)$ is
    \emph{stationary over $A$}. That is, there exists a unique
    global type extending it which does not fork over $A$. In the
    example above, the first global type is determined by a
    decreasing Morley sequence (which is a co-heir sequence), and
    the second one - by an increasing one.

    Having realized that $\Av(I,A \cup I)$ has a unique extension
    $q$ which does not fork over $A$, we proceed to understanding $q$.
    The natural conjecture that $q$ is the global average
    type of $I$ (as was the case in \cite{Us}) fails immediately:
    e.g., taking $I$ to be an increasing sequence in the type ``at
    infinity'' over $(\setQ,<)$ discussed above. Instead, we have to
    work with the so-called ``eventual'' type of $I$, $\Ev(I)$. This
    notion is essentially due to Poizat \cite{Poi}, although the
    name was proposed by Adler in \cite{Ad}. Poizat studied eventual
    types of ``special'' sequences. We work with a slight
    generalization of his notion, which we call ''Lascar special'' or ``weakly special''
    sequences.

    Given a special sequence $I$ over a set $A$ (see definitions in section 2), Poizat gave a natural construction of
    a global type $\Ev(I)$ which extends $\Av(I,A\cup A)$ and does
    not split over $A$. Using a similar construction, given a weakly
    special sequence $I$ over $A$, we obtain a global type $\Ev(I)$
    which extends the average type of $I$ and does not \emph{Lascar
    split} (equivalently, does not \emph{fork}) over $A$. Then we
    show that a sequence $I$ is weakly special over $A$ if and only
    if it is Morley over $A$. Summarizing all of the above, we can conclude that:

    \begin{itemize}
    \item
        A Morley sequence over $A$ determines a unique global type
        which does not fork ( = is Lascar invariant) over $A$. This
        type is precisely the eventual type of $I$, where $I$ is viewed as a
        weakly special sequence over $A$.
    \end{itemize}

    The reader might ask at this point how the results in \cite{Us} fall into this picture. From the construction of
    the eventual type, it is
    easy to see that if $I$ is a Morley indiscernible set, then
    $\Ev(I) = \Av(I,\fC)$. So it follows that $\Av(I,\fC)$ does not
    fork over $A$. We go further and show that this
    characterizes indiscernible sets: $\Av(I,\fC)$ does not fork
    over $A$ iff $\Av(I,\fC) = \Ev(I)$ iff $I$ is an indiscernible
    set iff $\Av(I,A)$ is generically stable. This gives us another
    nice characterization of generically stable types.

    Although all Morley sequences are important objects, some of
    them might be more useful than others. Recall that one of the
    most important properties of Morley sequences in simple theories
    is the following fact (which we refer to as ``Kim's Lemma''):
    Suppose the formula $\ph(x,a)$ divides over a set $A$; then for
    every Morley sequence $I$ in $\tp(a/A)$, we have that $\ph(x,I)
    = \set{\ph(x,a')\colon a'\in I}$ is inconsistent. A natural
    question is: is some version of Kim's Lemma true in dependent
    theories? In \cite{OnUs1} Alf Onshuus and the author showed (roughly) that
    it becomes true if one replaces ``for all'' Morley sequence with
    ``there exists'' one. The proof works, in fact, in a more
    general context of $NTP_2$ theories which includes both
    dependent and simple theories. Still, for many applications this
    result is insufficient. So one could ask: is there a stricter
    notion of a Morley sequence for which Kim's Lemma is true is
    stated?

    We answer this question positively in section 4 of the
    paper, using the notion of ``strict nonforking'' introduced by
    Shelah in \cite{Sh783}. That is, we show that for strict Morley
    sequences Kim's Lemma is true. We also investigate the notion of
    strict nonforking and show that a sequence in a global
    nonforking heir is a strict Morley sequence. Furthermore, we
    show existence of such sequences over models by deducing existence of
    nonforking heirs from recent results of Chernikov and Kaplan.
    Let  us also point out that Kim's Lemma for strict Morley
    sequences was shown by Chernikov and Kaplan independently in
    \cite{ChKa}.

    We conclude the paper with further properties of strict nonforking. In
    particular we show  weak version of local character (``bounded
    weight'' for strictly nonforking sequences) and discuss
    different versions of weak orthogonality.

    \medskip

    \emph{A note on terminology:} in \cite{Us} we restricted the
    term ``Morley sequence'' to a sequence in a definable type
    constructed with respect to a definition. Since then it has
    become very clear that arbitrary nonforking sequences are
    central objects, hence deserve a special name. Any nonforking (indiscernible) sequence
    (not necessarily in a definable type) is called ``Morley'' in a simple
    theory, so we decided to adopt the name. Note that in \cite{Us} we
    only worked with indiscernible sets, and it was shown there that
    any such nonforking sequence is also a sequence with respect to
    a definition. Hence the terminology here does not differ, in fact,
    from the one in \cite{Us}.

    \medskip
    \emph{A note on the framework:} Several proofs in this paper are
    carried out in the context of continuous model theory. We do
    this for several reasons. First, it has recently become clear
    that continuous model theory might be quite useful for studying
    dependent theories. Second, the terminology and the notation of
    continuous model theory are often very convenient. Third, some
    concepts developed in this paper require a slightly more
    sophisticated treatment when working in the continuous context,
    and we prefer to make sure everything works in this generality.

    In spite of all this, the main examples, motivations, etc,
    behind the results in this paper come from classical model
    theory. Hence the
    reader can safely assume (having gotten used to the slightly different notation, explained in the next subsection),
    that everything is happening in the
    classical (discrete) context. In fact, most of the proofs would
    not change at all if we decided to eliminate all traces of
    continuous model theory from the paper. The only subsection
    which would be significantly simplified is the discussion of
    the eventual type, subsection 2.2. Still, since we believe that eventual type is a central
    concept, we decided to develop it in the slightly higher generality.

\subsection{Notations}
In this paper, $T$ will denote a complete theory (sometimes
continuous), $\tau$ will denote the vocabulary of $T$, $L$ will
denote the language of $T$. We will assume that everything is
happening in the monster model of $T$ which will be denoted by
$\fC$. Elements and finite tuples of $\fC$ will be denoted $a,b,c$,
sets (which are all subsets of $\FC$) will be denoted $A,B,C$, and
models of $T$ (which are all elementary submodels of $\FC$) will be
denoted by $M,N$, etc.

Given an order type $O$ and a sequence \inseq{a}{i}{O}, we often
denote $a_{<i} = \lseq{a}{j}{i}$, similarly for $a_{\le i}$,
$a_{>i}$, etc. We will often identify a tuple $a$ or a sequence
$\inseq{a}{i}{O}$ with the set which is its union, but it should
always be clear from the context what we mean (although sometimes
when confusions might arise, we make the distinction, e.g.
$\Av(I,\cup I)$ will denote the average type of a sequence $I$ over
itself).

By $a \equiv_A b$ we mean $\tp(a/A) = \tp(b/A)$. By $a
\equiv_{\lstp,A} b$ we mean that $a$ and $b$ have the same Lascar
strong type over $A$; we also write $\lstp(a/A) = \lstp(b/A)$. By $a
\equiv_{ind,A} b$ we mean that there is an $A$-indiscernible
sequence containing both $a$ and $b$. Obviously $a \equiv_{ind,A}b
\Longrightarrow a \equiv_{\lstp,A}b \Longrightarrow a\equiv_Ab$.

Since part of the time we are working with a continuous theory, we
adopt some notation that seems convenient. For example, we denote
the truth value of a formula $\ph(x,a)$ with respect to a type $p
\in S(A)$ with $a \in A$ by $\ph^p(x,a)$. When no confusion arises,
a sentence $\psi(a)$ will denote its truth value in the monster
model. So $\ph^p(x,a) = \ps(b)$ means, when working with a classical
first order theory, that for every (some) realization $c \models p$
we have $\fC\models \ph(c,a) \leftrightarrow \ps(b)$. Given a
sequence of truth values $\lseq{t}{i}{\lam}$ which is eventually
constant, we denote the value which appears co-boundedly many times
by $\lim_{i<\lam}t_i$ (clearly, this is the limit of the sequence in
discrete topology).

\subsection{Preliminaries}

Recall that a theory (discrete or continuous) $T$ is called
\emph{dependent} if for every
    indiscernible sequence $I = \lseq{a}{i}{\lam}$, a formula
    $\ph(x,y)$ and $c$ we have
    $$\exists\lim_{i<\lam}\ph(a_i,c)$$

Recall that a sequence $I = \inseq{a}{i}{O}$ (where $O$ is a linear
ordering) is called \emph{indiscernible} over a set $A$ if the type
of $a_{i_1},\ldots,a_{i_k}$ over $A$ depends only on the order
between the indices $i_1, \ldots, i_k$ for every $k$. $I$ is called
an indiscernible \emph{set} if the type above depends on $k$ only.

We will call an \om-type $Q(x_0,x_1,\ldots)$ \emph{indiscernible}
over a set $A$ if every (equivalently, some) realization of it is an
indiscernible sequence over $A$. Note that $x_i$ can be finite
tuples. Clearly by compactness we can speak (slightly abusing the
usual terminology) of realizations $I$ of $Q$ where $I$ has any
(infinite) order type. In other words, we treat $Q$ as the
Ehrenfeucht-Mostowski type of an indiscernible sequence. We will
always assume, though, that our indiscernible sequences \emph{do not
have a last element}.

Slightly abusing notation, given an $A$-indiscernible sequence $I$,
we will often say that $Q$ as above is the type of $I$ over $A$ even
if $I$ is not of order type $\om$; similarly for Lascar type. Given
two such sequences $I$ and $J$, we write $I \equiv_A J$ if they are
realizations of the same indiscernible (over $A$) type $Q$ (but
possibly not of the same order type); similarly for Lascar type. In
other words, we write $I \equiv_A J$ for $\EM(I/A) = \EM(J/A)$.

\begin{dfn}
    Let $I$ be an $A$-indiscernible sequence. We say that a sequence
    $J$ \emph{continues} $I$ over $A$ if $I^\frown J$ is
    $A$-indiscernible.
\end{dfn}

A hyperimaginary element $a$ is said to be \emph{bounded} over a set
$A$ if the orbit of $a$ under the action of $\Aut(\fC,A)$ is small,
i.e. of cardinality less than $|\fC|$. The \emph{bounded closure of}
$A$, denoted by $\bdd(A)$, is the collection of all elements bounded
over $A$.


\subsection{Global Assumptions}

    All theories mentioned here are \emph{assumed to be
    dependent} unless stated otherwise. For the sake of
    clarity of presentation we also assume $T = T^{eq}$.


\subsection{Forking and splitting in dependent theories}

The following observations are well-known by now. The proofs can be
found e.g. in \cite{Us}, section 2.

\begin{fct}\label{fct:strongsplit}
\begin{enumerate}
\item
    Strong splitting implies dividing.
\item
    Lascar splitting implies forking.
\item
There are boundedly many global types which do not fork over a given
set $A$.
\end{enumerate}
\end{fct}

\begin{fct}\label{fct:forkind}
\item
    Let $I = \lseq{a}{i}{\lam}$ be such that
    \begin{itemize}
    \item
        $\tp(a_i/Aa_{<i})$ does not fork over $A$
    \item
        $\lstp(a_i/Aa_{<i}) = \lstp(a_j/Aa_{<i})$ for every $j \ge
        i $.
    \end{itemize}
    Then $I$ is a Morley (nonforking) sequence over $A$ (that is, it is
    indiscernible over $A$).
\end{fct}


\begin{fct}\label{fct:globalfork}
The following are equivalent for a global type $p$ and a set $A$:
\begin{itemize}
\item
    $p$ forks over $A$
\item
    $p$ divides over $A$
\item
    $p$ splits strongly over $A$
\item
    $p$ Lascar splits over $A$
\end{itemize}
\end{fct}

We will use the facts above all the time, sometimes without quoting.



\section{Average and Eventual types of indiscernible sequences}


\subsection{Average types}

Proofs in this subsection are carried out in the setting of
continuous logic, but the reader can easily ignore this and think
only in terms of classical model theory.

\begin{dfn}
    Let $I = \lseq{a}{i}{\lam}$ be an indiscernible
    sequence
    $B$ a set. We define the \emph{average type}
    of $I$
    over $B$ to be
    $$\Av(I,B)  = \set{\ph(x,b) = \lim_i\ph(a_i,b)}$$.
\end{dfn}

Recall

\begin{fct}
    $T$ is dependent if and only if for every $I$, $B$ as above, we have $\Av(I,B) \in \tS(B)$.
\end{fct}

\begin{rmk}\label{rmk:avind}
    Let $I$ be an indiscernible sequence over a set $A$. Then $a \models \Av(I,A\cup I)$
    if and only if $I^\frown\set{a}$ is indiscernible over $A$.
\end{rmk}

\begin{obs}\label{obs:concat_indisc}
    Let $I = \inseq{a}{i}{O}$ be an indiscernible sequence over a set $A$ and let $p$
    be a global type which extends $\Av(I,A\cup I)$ and does not fork
    over $A$. Suppose that $I' = \inseq{a'}{i}{O'}$ satisfies
    $a'_i \models p\rest AIa'_{<i}$. Then $J = I^\frown I'$ is
    indiscernible over $A$.
\end{obs}
\begin{prf}
    The mere existence of $p$ implies that $I$ is a nonforking
    sequence.
    Clearly $I'$ is
    nonforking over $AI$ based on $A$. Hence by Fact \ref{fct:forkind} it is enough to show that
    $a'_j \equiv_{\lstp,a_{<i}} a_i$ for all $i,j<\om$. But this is
    also clear since $a'_j \models \Av(I,Aa_{<i})$, so
    $I^\frown{a'_j}$ is indiscernible, hence in fact $a'_j \equiv_{ind,a_{<i}}
    a_i$.
\end{prf}

\begin{obs}\label{obs:concat_element}
    Let $I = \inseq{a}{i}{O}$ be an indiscernible sequence over a set $A$ and let $p$
    be a global type which extends $\Av(I,A\cup I)$ and does not fork
    over $A$. Suppose that $I' \equiv_{\lstp,A} I$. Then $p \rest
    AI' = \Av(I',A \cup I')$.
\end{obs}
\begin{prf}
    As before, existence of $p$ implies that $I$ is a nonforking
    sequence. Clearly $I'$ is
    nonforking over $A$.

    Let $c' \models p\rest AI'$; clearly, it is enough to show that
    $I'^\frown{c'}$ is $A$-indiscernible. Let $c$ be such that
    $Ic \equiv_{\lstp,A} I'c'$; it is enough to show that $I^\frown{c}$
    is $A$-indiscernible. Note that $p$ does not fork, hence does
    not Lascar-split over $A$. It follows that $c \models p\rest
    AI = \Av(I,A\cup I)$, and we are done.
\end{prf}

\begin{rmk}\label{rmk:weaklyspecial}
    Note that in the previous Observation, $p$ extends both the averages of $I$ and $I'$ over $A$.
    Moreover, if in the proof of the Observation we chose $c \models p\rest AII'$, we would have
    that both $I^\frown{c}$ and $I'^\frown{c}$ are
    $A$-indiscernible. This property of a Morley sequence will become very important
    later.
\end{rmk}

\begin{obs}\label{obs:avextend}
\begin{enumerate}
\item
    Let $I$ be an indiscernible sequence over a set $A$, $p$ a global
    type extending $\Av(I,A\cup I)$ which does not fork over $A$.
    Then for every $A$-indiscernible sequence $I'$ continuing $I$, we
    have $p\rest AII' = \Av(I',AII')$.
\item
    Same conclusion if $p$ is just a type over $AII'$ which does not
    split strongly over $A$.
\end{enumerate}
\end{obs}
\begin{prf}
\begin{enumerate}
\item
    Let $I = \inseq{a}{i}{O}$, $I' = \inseq{a'}{i}{O'}$.
    Assume that $\ph(x,a_{<i}a'_{<i}) \in \Av(J,AIJ)$ where
    $\ph(x,yy')$ is a formula over $A$. Then clearly since $I'$
    continues $I$ it is also the case that $\ph(x,a_{<2i}) \in
    \Av(I,AI) \subseteq p$. As $p$ does not fork, hence does not split strongly over $A$,
    we have $\ph^p(x,a_{<i},a'_{<i}) = \ph^p(x,a_{<2i})$, as
    required.
\item
    Same proof.
\end{enumerate}
\end{prf}

We have mentioned that a type $p \in \tS(A)$ in a dependent theory
has boundedly many global nonforking extensions. This does not mean,
of course, that $p$ is stationary (even if $A$ is a model), unlike
in the stable case. The following lemma shows that once given a
Morley (nonforking) sequence in $p$, it completely determines a
global nonforking extension.

\begin{lem}\label{lem:forkseqstat}
\begin{enumerate}
\item
    Let $I$ be an indiscernible sequence over a set $A$ and
    let $p,q$ be global types extending $\Av(I,A\cup I)$, both
    do not fork over $A$. Then $p=q$.
\item
    Let $I$ be a Morley (nonforking) sequence over a set $A$. Then there exists a unique
    global types extending $\Av(I,A\cup I)$ which does not fork over
    $A$. In other words, $\Av(I,A\cup I)$ is stationary over $A$.
\end{enumerate}
\end{lem}
\begin{prf}
\begin{enumerate}
\item
    Assume towards contradiction that $q \neq
    p$, so there is $\ph(x,b)$ such that $\ph(x,b)\in p$ but
    $\neg\ph(x,b)\in q$.

    Construct by induction on $\al<\om$ sequence $J_\al = \lseq{a^\al}{i}{\om}$ such that
    \begin{itemize}
    \item
        $a^{2\al}_i \models p \rest Ab IJ_{<\al}a^{2\al}_{<i}$
    \item
        $a^{2\al+1}_i \models q \rest Ab IJ_{<\al}a^{2\al+1}_{<i}$
    \end{itemize}

    We claim that $J = J_0^\frown J_1^\frown\cdots $ is an indiscernible sequence. Once
    we have shown this, it yields an immediate contradiction to dependence, since
    $\ph(a^\al_i,b) \neq  \ph(a^{\al+1}_i,b)$ for all $\al$.

    So we show by induction on $\al$ that $J^\al = I^\frown J_0^\frown
    \cdots^\frown J_\al$ is indiscernible (even over $A$). For $\al
    = 0$ this is true by Observation \ref{obs:concat_indisc}.

    Let us take care of $\al = 1$ (the continuation is the same). By
    Observation \ref{obs:avextend} $q$ extends $\Av(J^0,A\cup J^0)$.
    Now apply Observation \ref{obs:concat_indisc} again.

\item
    This is just a restatement of (i).
\end{enumerate}
\end{prf}

\begin{exm}
    Let $p$ be the type ``at infinity'' over $(\setQ,<)$. Then it
    has \emph{two} global nonforking extensions: one is the type of
    the cut $\setQ^+$, which is finitely satisfiable in $\setQ$, and
    the other one is the type ``at infinity'' over the monster
    model, which is the extension with respect to the definition of
    $p$. The first one is determined by a decreasing Morley
    sequence, whereas the second one - by an increasing one (those
    are precisely the two possible types of a Morley sequence in
    $p$).
\end{exm}

\subsection{Eventual types}

So we have shown that a nonforking sequence $I$ over $A$, it
determines uniquely a global type that does not fork over $A$. It is
natural to ask what this global extension looks like. A first guess
might be $\Av(I,\fC)$ (after all, this is the case in stable
theories), but it is easy to see that this does not always work. For
example, taking $I$ to be an increasing sequence of rational numbers
in the structure $(\setQ,<)$, which is clearly nonforking over the
empty set, the global nonforking extension is the type ``at
infinity'' and not the global average (which is the type of the cut
of $I$).

The reader might recall from section 2 of the authors previous
article \cite{Us} that when one starts with a ``stable-like'' (more
precisely, generically stable) type (so when $I$ is an indiscernible
\emph{set}) then $\Av(I,\fC)$ is indeed the unique global nonforking
extension of $\Av(I,A\cup I)$. But this case is in a sense ``too
easy'' and does not reflect the subtleties of the general situation.
As a matter of fact, we will see later in the paper that the global
average is the unique nonforking extension \emph{if and only if} it
is generically stable if and only if $I$ is an indiscernible set.

So for the general case we will need to apply a slightly more
careful analysis and introduce a different notion of a ``limit''
type of a nonforking sequence. The ideas behind this notion
generalize Poizat \cite{Poi}, Shelah \cite{Sh783} and, more
recently, Adler \cite{Ad}, who work with co-heir and ``special''
(see below) sequences, but as we shall see, everything generalizes
to an arbitrary Morley sequence quite easily.
To the best of our knowledge, Adler was the first to actually give
the central notion of this subsection a name. We will follow his
terminology and call the ``limit'' type we are interested in ``the
eventual type of the sequence''. In this subsection we work in the
continuous setting. This is the only place in the paper where the
proofs would become somewhat simpler if we decided to remain in the
classical context. But since the concept of eventual type seems
important enough, we decided to make an effort and do everything
carefully in the more general framework.


\begin{obs}\label{obs:alternate}
  For every $\ph(x,y)$ (maybe with parameters, maybe $y$ is empty)
  and $\eps>0$
  there exists $k <\om$ such that
  there does not exist an infinite indiscernible sequence
    $\lseq{b}{i}{\om}$ and
    $c \in \fC$,
    such that for all $i<k$ we have
    $$ |\ph(b_{i},c)-\ph(b_{i+1},c)|> \eps$$
    We denote \emph{minimal} such $k$ by $k_{\ph,\eps}$.
\end{obs}

For notational convenience, let us introduce the following notation
for the number for $\eps$-alternations (see also Adler \cite{Ad} for
a related notion of ``alternation rank'').

\begin{dfn}
\begin{enumerate}
\item
    Let $I = \inseq{a}{i}{O}$ be an indiscernible sequence, $\ph(x,b)$ a formula,
    $\eps>0$. We denote by $\alt(\ph(x,b),\eps,I)$ the maximal $k$
    such that there exist $i_0 < \ldots<i_k \in O$ such that
    $|\ph(a_{i_j},b)-\ph(a_{i_{j+1}}, b)|>\eps$ for all $j<k$.
    Clearly $\alt(\ph(x,b),\eps,I) \le k_{\ph(x,b),\eps}$.
\item
    When we omit $I$, we mean the maximum over all $I$,
    that is, $\alt(\ph(x,b),\eps) = k_{\ph(x,b),\eps}$.
\item
    When we replace $I$ with an indiscernible type $Q$, we mean
    the maximum over all realizations of $Q$.
\end{enumerate}
\end{dfn}

The following notion is due to Poizat \cite{Poi}.

\begin{dfn}
\begin{enumerate}
\item
    We call an $A$-indiscernible type sequence $I$ \emph{special} if for every
    two realizations $I_1$ and $I_2$ of $\tp(I/A)$, there exists $c$ such that
    $I_1^\frown c$ and $I_2^\frown c$ are $A$-indiscernible.
\item
    In this case we call the type $\tp(I/A)$ \emph{special} over
    $A$.
\end{enumerate}
\end{dfn}

It is known in classical model theory (see e.g. \cite{Ad}) that a
sequence $I = \inseq{a}{i}{O}$  is special over $A$ if and only if
there exists a global type $p$ which does not split over $A$ and
$a_i \models p\rest Aa_{<i}$ for all $i \in O$.
The same proof works in the continuous context (as we shall see),
but this is not quite what we are looking for: the assumption that
$p$ does not split over $A$ is too strong for us, and we would like
to replace it with forking, equivalently, Lascar splitting. This is
why we will find it more convenient to work with a Lascar strong
type of a sequence instead of a type.

\begin{dfn}
\begin{enumerate}
\item
    We call an $A$-indiscernible sequence \emph{weakly special} (or \emph{Lascar special}) if
    for any two realizations $I_1$ and $I_2$ of $\lstp(I/A)$, there exists $c$ such that
    $I_1^\frown c$ and $I_2^\frown c$ are $A$-indiscernible.
\item
    In this case we say that the type $\tp(I/A)$ is \emph{weakly special} (or \emph{Lascar special}) over
    $A$ and that Lascar strong type $\lstp(I/A)$ is \emph{special}
    over $A$.
\end{enumerate}
\end{dfn}

The following follows by induction and compactness:

\begin{obs}
    Let $Q$ be a Lascar special type over $A$. Then for every collection
    $\lseq{I}{i}{\lam}$ of $I_i \models Q$ of the same Lascar strong type over $A$
    and for every order type $O$ there exists
    $J \models Q$ of order type $O$ which continues all the $I_i$'s over
    $ A$.
\end{obs}

\begin{obs}\label{obs:max_special}
    Let $Q$ be a Lascar special type over $A$, $\ph(x,b)$ a formula,
    $\eps>0$, $J,J' \models Q$, $J \equiv_{\lstp,A}J'$, such that $\alt(\ph(x,b),\eps,J) =
    \alt(\ph(x,b),\eps,J') = \alt(\ph(x,b),\eps,Q)$. Denote
    $q=\Av(J,Ab)$, $q' = \Av(J',Ab)$. Then
    $$ |\ph^q(x,b)-\ph^{q'}(x,b)|\le 2\eps$$
\end{obs}
\begin{prf}
    Denote $\rho = \ph^q(x,b)$, $\rho' = \ph^{q'}(x,b)$. So
    $\lim_J\ph(x,b) = \rho$, $\lim_{J'}\ph(x,b) = \rho'$.
    Let $a$ be such that both
    $J^\frown a$ and $J'^\frown a$ are indiscernible. Clearly both
    satisfy $Q$, and so by the assumption the value
    of $\ph(x,b)$ can not change by more than $\eps$ in any
    sequences extending $J$ or $J'$. So $|\rho-\ph(a,b)|\le\eps$
    and $|\rho'-\ph(a,b)|\le\eps$, hence $|\rho-\rho'|\le 2\eps$
    as required.
\end{prf}

\begin{obs}
    Let $I$ be a Lascar special sequence over $A$, $\ph(x,b)$
    a formula, $\eps>0$. Let $J \equiv_{\lstp,A}I$ be such that
    $\alt(\ph(x,b),\eps,J) = \alt(\ph(x,b),\eps,Q)$. Then $I$
    can be extended to $I'$ such that $|\lim_{I'}\ph(x,b) -
    \lim_J\ph(x,b)|\le\eps$.
\end{obs}
\begin{prf}
    Let $J'\models \lstp(I/A)$ be an \om-sequence which
    continues both $I$ and $J$; clearly,
    $$|\lim_{J'}\ph(x,b)-\lim_J\ph(x,b)|\le\eps$$.
\end{prf}

\begin{dfn}
    Let $\ph(x,b)$ be a formula. We say that an indiscernible
    sequence $J$ \emph{eventually determines} $\ph(x,b)$ if
    $\lim_{J'}\ph(x,b)$ is constant for all $J'$ continuing $J$.
\end{dfn}

\begin{lem}\label{lem:evtype}Let $Q$ be a Lascar special type over $A$, $\ph(x,b)$ a formula.
\begin{enumerate}
    \item
    There exists $J \models Q$ which eventually determines
    $\ph(x,b)$.
\item
    For every $I,J\models Q$ with $J \equiv_{\lstp,A}I$ which eventually determine $\ph(x,b)$
    we have $\lim_I\ph(x,b) = \lim_J\ph(x,b)$, that is, the
    ``eventual value'' of $\ph(x,b)$ depends only on
    Lascar strong type of $J$ over $A$, and not
    on the choice of $J$. We call this number the \emph{eventual
    value of $\ph(x,b)$ with respect to $\lstp(J/A)$}.
\item
    Every $I \models Q$ can be extended to $J \models Q$ that
    eventually determines $\ph(x,b)$. For every $J,J'$ which
    continue $I$ and eventually determine $\ph(x,b)$
    we have $\lim_J\ph(x,b) = \lim_J'\ph(x,b)$.
\end{enumerate}
\end{lem}
\begin{prf}
\begin{enumerate}
\item
    Let $\eps>0$, and let $J_\eps\models Q$ be such that
    $\alt(\ph(x,b),\eps,J) = \alt(\ph(x,b),\eps,Q)$ and denote
    $\rho_\eps = \lim_{J_\eps}\ph(x,b)$.

    Let $\zeta<\eps$. As in previous Observation we can extend $J_\eps$ to $J'_{\eps}$ such that
    $|\lim_{J'_\eps}\ph(x,b)-\lim_{J_\zeta}\ph(x,b)|\le\zeta$.
    Also, clearly
    $|\lim_{J'_\eps}\ph(x,b)-\lim_{J_\eps}\ph(x,b)|\le\eps$. So
    $|\rho_\zeta-\rho_\eps|\le\zeta+\eps$.

    Denote $J_n = J_{\frac 1n}$, $\rho_n = \rho_{\frac 1n}$. The
    inequality above amounts to $|\rho_n-\rho_m|\le \frac 1n+\frac
    1m$. So $\rho_n$ converges; denote $\rho = \lim \rho_n$.

    Note that we may assume that if $\zeta<\eps$ then $J_\zeta$
    extends $J_\eps$. This is because we can replace $J_\zeta$ with
    $J_\eps J'$ where $J'$ is an \om-sequence which continues both $J_\zeta$ and $J_\eps$.
    So working only with $\eps_n = \frac 1n$, we can do
    the replacement process above for all $n$ by simple induction
    on the natural numbers, obtaining $J_n  \subseteq
    J_{n+1}$, and for simplicity assume $\otp(J_n) =
    n\om$.

    By compactness there is $J$ of order type $\om^2$
    which extends all the $J_\eps$'s. Clearly, if $J'$ extends $J$,
    then the value of $\lim_{J'} \ph(x,b)$ has to agree with
    $\lim_{J_n}\ph(x,b) = \rho_{\frac 1n}$ up to $\frac 1n$,
    hence has to equal $\rho$, which completes the proof of the first part of the
    claim.
\item
    Now suppose $J$ and $I$ both eventually determine $\ph(x,b)$
    and satisfy the same Lascar strong type. We can find $J'$ which continues both; so
    $$ \lim_I \ph(x,b) = \lim_{IJ'}\ph(x,b) =
    \lim_{J'}\ph(x,b) = \lim_{JJ'}\ph(x,b) = \lim_J\ph(x,b)$$
    as required.
\item
    Let $J' \models Q$ eventually determine $\ph(x,b)$ (see clause
    (i) above), and let $J$ continue both $I$ and $J'$, clearly $J$
    is as required.

    For uniqueness, note that $\lstp(J/A) = \lstp(I/A)$ for every
    $J$ continuing $I$, and use clause (ii).
\end{enumerate}
\end{prf}

We are now ready to make the main definition of the subsection.

\begin{dfn}\label{dfn:evtype}
\begin{enumerate}
\item
    Given a Lascar special sequence $I$ over $A$ and a formula $\ph(x,b)$ we
    denote the \emph{eventual value} of $\ph(x,b)$ with respect to
    $I$ (see Lemma \ref{lem:evtype} (ii)) by $\ph^I(x,b)$. If $Q$ is a type which implies $\lstp(I/A)$ (e.g.
    $Q = \tp(I/M)$ where $M$ is a model containing $A$, or just
    $Q = \tp(I/\bdd^{heq}(A))$), we write $\ph^Q(x,b)$ for
    $\ph^I(x,b)$ (clearly, $\ph^I(x,b)$ depends only on $Q$).
\item
    Given a Lascar special sequence $I$ over $A$ and a set $C$ we define the \emph{eventual type} of
    $I$ over $C$, $q=\Ev(I,C)$ as follows: given a formula
    $\ph(x,b)$ over $C$, let $\ph^q(x,b)$ equal $\ph^I(x,b)$.

    Again, if $Q$ implies $\lstp(I/A)$, we write $\Ev(Q,C)$.
\item
    If we omit $C$, we mean $C = \fC$ (so we obtain a global type).
\end{enumerate}
\end{dfn}

\begin{rmk}
    $\Ev(I,C)$ as above is well-defined and is a complete type over $C$.
\end{rmk}
\begin{prf}
    Compactness and existence + uniqueness of eventual value, that is, Lemma
    \ref{lem:evtype}.
\end{prf}

\begin{rmk}\label{rmk:evav}
    Let $I$ be a Lascar special sequence over $A$. Then
    $\Ev(I,A\cup I) = \Ev(I)\rest AI = \Av(I,A\cup I)$.
\end{rmk}
\begin{prf}
    So $a \models \Ev(I,A\cup I)$ if and only if for every formula $\ph(x,a_{<i})$ over $A$ and for
    some/every
    continuation $J$ of $I$ which eventually determines
    $\ph(x,a_{<i})$ we have
    $$\ph(a,a_{<i}) = \ph^J(x,a_{<i}) = \lim_J\ph(x,a_{<i}) = \lim_I\ph(x,a_{<i})$$
    The last equality is true because $J$ is an indiscernible
    sequence continuing $I$. So clearly $a \models\Ev(I,A\cup I)$
    if and only if $a\models\Av(I,A\cup I)$.
\end{prf}

\begin{obs}\label{obs:evnonsplit}
    Let $I$ be Lascar special over $A$. Then $\Ev(I)$ does not Lascar-split
    (equivalently,
    does not fork) over $A$.
\end{obs}
\begin{prf}
    Let $b \equiv_{\lstp,A} b'$, $\ph(x,y)$ a formula. Suppose that the
    eventual value of $\ph(x,b)$ is determined by $J \equiv_{\lstp,A} I$.
    Let $J'$ be such that $J'b' \equiv_{\lstp,A} Jb$. Clearly $J'$
    eventually determines $\ph(x,b')$, $J' \models Q$ and
    $\lim_J\ph(x,b) = \lim_{J'}\ph(x,b')$. By uniqueness of
    eventual value, $J'$ determines
    $\ph^I(x,b')$, so we are done.
\end{prf}

\begin{dsc}
    Note that we could have gone through the same process with
    special sequences instead of Lascar special, replacing ``Lascar
    strong type'' with ``type'' in all the statements and proofs
    above. In this case, the eventual type of a special sequence $I$
    over $A$ depends only on $Q = \tp(I/A)$, and we denote it by $\Ev(Q)$. In the Observation
    above we would get then:
\end{dsc}

\begin{obs}
    Let $Q$ be a special type over $A$. Then $\Ev(Q)$ does not split (and therefore
    does not fork) over $A$.
\end{obs}

\begin{cor}\label{cor:evnonsplit}
\begin{enumerate}
\item
    Let $I$ be Lascar special over $A$. Then for every set $C$ the type
    $\Ev(I,C)$ does not fork over $A$.
\item
    Let $Q$ be special over $A$. Then for every set $C$ the type
    $\Ev(Q,C)$ does not either split or fork over $A$.
\end{enumerate}
\end{cor}

The following is a known characterization of special sequences
mentioned above:

\begin{lem}
    An $A$-indiscernible type $Q(x_i)$ is special if and only if there
    exists a global type $q = q(x)$ which does not split over $A$
    such that any realization $I = \inseq{a}{i}{O}$
    of $Q$ satisfies $a_i \models q\rest Aa_{<i}$.
\end{lem}
\begin{prf}
    If any realization of $Q$ is a nonsplitting sequence in $q$ over
    $A$, then given $I,I'\models Q$ simply choose $c\models q\rest
    AII'$. By nonsplitting both $I^\frown c$ and $I'^\frown c$ are
    indiscernible.

    On the other hand, assume that $Q$ is special over $A$ and let
    $q = \Ev(Q)$. Then by Corollary \ref{cor:evnonsplit} $q$ does not split over $A$ and clearly for any
    realization $I$ of $Q$ and for every $a_i \in I$ we have $a_i
    \models \Ev(Q)\rest Aa_{<i}$.
\end{prf}

We are more interested in the similar characterization of Lascar
special sequences:

\begin{lem}
    An $A$-indiscernible sequence $I$  is Lascar special if and only if
    it is a Morley (nonforking) sequence over $A$.
\end{lem}
\begin{prf}
    A nonforking sequence is Lascar special by Remark
    \ref{rmk:weaklyspecial}.

    On the other hand, if $I$ is Lascar special, then by Observation \ref{obs:evnonsplit} $q = \Ev(I)$
    does not fork over $A$, hence $\Av(I,A\cup I) = \Ev(I,A\cup I)$
    (see Remark \ref{rmk:evav}) does not fork over $A$; so $I$ is nonforking.
\end{prf}

\begin{dsc}
    Comparing the two lemmas above, the reader can notice that
    unlike special sequences, which are characterized in terms of
    existence of a certain \emph{global} type, Lascar special types
    have an \emph{internal} characterization (it can be seen from
    the type whether or not it is nonforking). This is part of the
    reason we believe that Lascar special (equivalently, Morley)
    sequences are a better notion.

    Note that simply replacing ``nonforking'' with ``nonsplitting''
    in the characterization of Lascar special sequences will
    \emph{not} lead to a characterization of special sequences. This
    is because nonsplitting does not satisfy the extension axiom; in
    fact, locally nonsplitting does not imply nonforking (unless one
    works over slightly saturated models).
    See section 6 of \cite{Us} for examples of forking nonsplitting
    sequences.
\end{dsc}

\begin{prp}\label{prp:evstat}
\begin{enumerate}
\item
    Let $I$ be a Lascar special sequence over $A$. Then $\Ev(I)$ is the
    unique global extension of $\Av(I,  A\cup I)$ that does not split over
    $A$.
\item
    Let $Q$ be a special type over $A$, $I \models Q$. Then $\Ev(Q)$ is the
    unique global extension of $\Av(I,A\cup I)$ that does not split over
    $A$.
\end{enumerate}
\end{prp}
\begin{prf}
    We have already shown in Observation \ref{obs:evnonsplit}
    that $\Ev(Q)$ does not split over $A$. By
    Remark \ref{rmk:evav} $\Ev(Q)$ extends $\Av(I,A\cup I)$.
    Uniqueness follows now from Lemma \ref{lem:forkseqstat}.

\end{prf}

\section{Indiscernible sets and generic stability}

In this section we will characterize generically stable types in
terms of the eventual type of their Morley sequences. The proofs are
carried out in the continuous context.

\subsection{Indiscernible sets}

Having achieved in the previous section a pretty good understanding
of an arbitrary nonforking sequence $I$, we will try to draw more
conclusions assuming that $I$ is an indiscernible set.

Recall

\begin{dfn}
  Let $I = \inseq{b}{i}{O}$ be an infinite indiscernible sequence. We say that
  a formula $\ph(x,y)$ is \emph{stable} for $I$ if for every $c \in
  \fC$ and $\eps>0$, either
  the set $\set{i \in I \colon \ph(b_i,c)<\eps}$
  or the set $\set{i \in I \colon \ph(b_i,c)>\eps}$ is finite.
\end{dfn}

Restating Observation \ref{obs:alternate} for indiscernible sets we
get

\begin{obs}\label{obs:stabform}
  If $\inseq{b}{i}{I}$ is an infinite indiscernible set, then every $\ph(x,y)$
  is stable for $\seq{b_i}$. Moreover,
  for every $\ph(x,y)$ and $\eps>0$
  there exists $k = k_{\ph,\eps} < \om$ such that for every $c \in \fC$,
    for all but $k$-many $i<\lam$ we have
    $$ |\ph(b_{i_1},c)-\ph(b_{i_2},c)|\le \eps$$
\end{obs}

\begin{obs}
    Let $I_1,I_2$ be equivalent indiscernible sets, that is, there is
    an indiscernible sequence (set) $J$ continuing both.
    Then $\Av(I_1,\fC) = \Av(I_2,\fC)$.
\end{obs}
\begin{prf}
    Let $\ph(x,b)$ be a formula. By Observation \ref{obs:stabform}
    it is easy to see that
    $$\lim_{I_1}\ph(x,b) = \lim_{I_1J}\ph(x,b) =
    \lim_J\ph(x,b)$$
    and similarly for $I_2$.
\end{prf}

\begin{obs}
\begin{enumerate}
\item
    Let $Q$ be a special type of an indiscernible set. Then $\Ev(Q)
    = \Av(I,\fC)$ for any realization $I$ of $Q$.
\item
    Let $I$ be a nonforking sequence over $A$ which is also an
    indiscernible set. Then $\Ev(I) = \Av(I,\fC)$.
\end{enumerate}
\end{obs}
\begin{prf}
    For (i), let $I \models Q$, $\ph(x,b)$ a formula, $J$ continuing $I$
    eventually determines $\ph(x,b)$. By the previous observation
    $\lim_I\ph(x,b) = \lim_J\ph(x,b)$, as required.

    (ii) Similar.
\end{prf}

As a consequence we can conclude the main result of section 2 of
\cite{Us}:

\begin{cor}\label{cor:seqstat}
    Let $I = \lseq{b}{i}{\om}$ be a nonforking sequence over $A$ which is also an
    indiscernible set. Denote $p = \Av(I,A\union I)$. Then $p$ has
    a unique nonsplitting extension over $\fC$, which equals
    $\Av(I,\fC)$. In particular, $\Av(I,\fC)$ does not fork over
    $A$.
\end{cor}

\subsection{Generically stable types}

We recall the main equivalences and facts from \cite{Us}.


\begin{thm}
    The following are equivalent for a type $p \in \tS(A)$:
    \begin{enumerate}
    \item
        There is a nonforking sequence in $p$ which is an
        indiscernible set.
    \item
        $p$ is extensible and every nonforking sequence in it is an
        indiscernible set.
    \item
        $p$ is definable over $\acl(A)$ and some/every Morley
        sequence with respect to a definition of $p$ is an
        indiscernible set.
    \item
        Some/every global nonforking extension of $p$ is both definable
        over and finitely satisfiable in a model of density
        character
        $|A|+|T|$ containing $A$.
    \end{enumerate}
\end{thm}

In case one/all of the equivalences above hold/s, we call $p$
\emph{generically stable}.

\begin{thm}
    Let $p \in \tS(A)$ be a generically stable type definable over
    $A$ (e.g. $A = \acl(A)$, or $A$ is a Morley sequence in $p$, or
    just $p$ is finitely satisfiable in $A$, see section 5 of
    \cite{Us}). Then $p$ is stationary.
\end{thm}

We now connect eventual types to generic stability.

\begin{prp}\label{prp:evstab}
    Let $I$ be a Morley sequence over $A$, and assume that
    $\Ev(I) = \Av(I,\fC)$. Then $\Av(I,A)$ is
    generically stable, hence so is $\Av(I,\fC) = \Ev(I)$, and $I$ is an indiscernible set.
\end{prp}
\begin{prf}
    First, note that $q=\Av(I,\fC)$ does not split over $A$ (since it
    equals $\Ev(Q)$), so it is enough to show that $I$ is an
    indiscernible set. Without loss of generality we may assume
    $\otp(I) = \om$.

    Let $I = \lseq{a}{i}{\om}$. Construct $I' =
    \seq{a_{\om+i}\colon i<\om}$ by $a_i \models q\rest Aa_{<i}$.
    Clearly $II'$ is an $A$-indiscernible sequence, and it is enough
    to show that $I'$ is an indiscernible set.

    Assume $\ph(a_{<\om+i},a_{\om+i},a_{\om+i+1},a_{\om+i+2})=0$
    and let us show
    $\ph(a_{<\om+i},a_{\om+i+1},a_{\om+i},a_{\om+i+2})=0$. The
    general case will work in the same way. Let $\eps>0$.

    \begin{itemize}
    \item
        By indiscernibility
        $\ph(a_{<i},a_{\om+i},a_{\om+i+1},a_{\om+i+2})=0$.
    \item
        $\ph^q(a_{<i},a_{\om+i},a_{\om+i+1},x)=0$ (by the choice
        of $a_{\om+i+2}$)
    \item
        $q = \Av(I,\fC)$, so
        $\lim_I\ph(a_{<i},a_{\om+i},a_{\om+i+1},x)=0$.
    \item
        For all $k<\om$ big enough
        we have $\ph(a_{<i},a_{\om+i},a_{\om+i+1},a_k)\le \eps$.
    \item
        Similarly, $\ph(a_{<i},a_{\om+i},x,a_k) \in \Av(I,\fC)$,
        hence for all $k<\ell<\om$ big enough we have
        $\ph(a_{<i},a_{\om+i},a_\ell,a_k)\le 2\eps$.
    \item
        By indiscernibility
        $\ph(a_{<i},a_{\om+i+1},a_\ell,a_k)\le 2\eps$.
    \item
        By indiscernibility again (recall that $\ell>k$)
        we have $\ph(a_{<i},a_{\om+i+1},a_{\om+i},a_k)\le 2\eps$
        for all $k<\om$ big enough. So $\ph^q(a_{<i},a_{\om+i+1},a_{\om+i},x)\le
        2\eps$.
    \item
        $\ph^q(a_{<i},a_{\om+i+1},a_{\om+i},a_{\om+i+2})\le
        2\eps$.
    \end{itemize}
    Since $\eps$ was arbitrary, we are done.
\end{prf}

So we have obtained a new characterization of a generically stable
type. Note that clauses (iv) and (v) in the Theorem below with
``sequence'' replaced with ``set'' appears already in \cite{Us}. It
is interesting to find out that the requirement of $I$ being an
indiscernible set is not needed and follows from the fact that the
average type does not fork over $A$.

\begin{thm}
    Let $p \in \tS(A)$. Then the following are equivalent:
    \begin{enumerate}
    \item
        $p$ is
        generically stable.
    \item For some Morley sequence $I$ in $p$ we have $Ev(I) =
    \Av(I,\fC)$.
    \item $p$ is extensible and for all Morley sequences $I$ in $p$ we have $Ev(I) =
    \Av(I,\fC)$.
    \item For some indiscernible sequence $I$ in $p$, $\Av(I,\fC)$ does not fork over
    $A$.
    \item $p$ is extensible and for every Morley sequence $I$ in $p$, $\Av(I,\fC)$ does not fork over
    $A$.

    \end{enumerate}
\end{thm}
\begin{prf}
    \noindent (i) \then (iii) By Corollary \ref{cor:seqstat}.

    \noindent (iii) \then (ii) Clear.

    \noindent (ii) \then (i) By Proposition \ref{prp:evstab}.

    \noindent (ii) \then (iv) Clear by now.

    \noindent (iv) \then (ii) Clearly $\Av(I,\fC)$ extends $\Av(I,A\cup
    I)$; since it also does not fork over $A$, we immediately get
    that $I$ is a nonforking sequence, and by stationarity of
    $\Av(I,A\cup I)$ over $A$ (see Corollary \ref{cor:evnonsplit})
    $\Ev(I) = \Av(I,\fC)$.

    \noindent (iv) \iff (v) is easy by stationarity.
\end{prf}

\section{On strictly free extensions}

In this section we will investigate a strong notion of a Morley
sequence in a dependent theory, which turns out to have some very
good properties.

The following definitions are equivalent to the ones given in
section 5 of \cite{Sh783}.

\begin{dfn}
\begin{enumerate}
\item
    Let $A\subseteq B$. We say that a type $p \in \tS(B)$ is a strictly nondividing
    extension of $p\rest A$ if for every $a \models p$
    \begin{itemize}
    \item
        $\tp(a/B)$ does not divide over $A$
    \item
        $\tp(B/Aa)$ does not divide over $A$.
    \end{itemize}
\item
    Let $A\subseteq B$. We say that a type $p \in \tS(B)$ is a \emph{strictly
    free} (or strictly nonforking)
    extension of $p\rest A$ if there exists a global type $q$
    extending $p$ which is a strictly nondividing extension of
    $ p\rest A$. We also say that $p$ is \emph{strictly free} over
    $A$. If $a\models p$, we write $a\ind_A^{st}B$.
\end{enumerate}
\end{dfn}

\begin{rmk}
    We will say that a type $p \in \tS(B)$ \emph{co-divides} over
    a set $A$ if there is $a \models p$ such that $\tp(B/Aa)$
    divides  over $A$. In other words, $p$ co-divides over $A$ if
    there exist $a \models p$, a formula $\ph(x,b) \in p$ such that
    $\ph(a,y)$ divides over $A$.

    Clearly $p \in \tS(B)$ is a strictly non-dividing extension over $A$ if and
    only if it does not divide and does not co-divide over $A$. We
    will make several observations about co-dividing.
\end{rmk}

\begin{obs} (T any theory)
    Let $A \subseteq B \subseteq C$, $p \in \tS(C)$ is a heir of
    $p\rest B$ and $p \rest B$ does not co-divide over $A$. Then $p$ does not
    co-divide over $A$.
\end{obs}
\begin{prf}
    Assume towards contradiction that $a\models p$, $\ph(x,c) \in p$,
    $\ph(a,y)$ divides over $A$. Let $\ph(x,b) \in p\rest B$ (recall
    that $p$ is a heir over $B$). Clearly $a \models p\rest B$ and
    $\ph(a,y)$ exemplify co-dividing of $p \rest B$, a
    contradiction.
\end{prf}

\begin{obs} (T any theory)
    Let $A$ be a set, $N$ an $(|A|+|T|)^+$-saturated model
    containing $A$, $p \in \tS(N)$ does not split over $A$. Let $q$
    be the unique global extension of $p$ which does not split over
    $A$ (see e.g. \cite{Us}, Lemma 2.23). Then $q$ is a heir of $p$.
\end{obs}
\begin{prf}
    Let $\ph(x,c) \in q$, and choose $b \in N $ of the same type as
    $c$ over $A$. By nonsplitting, clearly $\ph(x,b) \in q\rest N =
    p$.
\end{prf}

\begin{cor}
    Let $N$ be saturated enough over a model $M$, $p \in \tS(N)$ is a strongly
    nondividing extension of $p \rest M$ (in particular $p$ does not split
    over $M$). Let $q$ be the unique global extension of $p$ which
    does not fork/split over $M$. Then $q$ is a strongly nondividing
    extension of $p \rest M$.
\end{cor}
\begin{prf}
    Clearly $q$ does not divide over $M$, so we only need to show
    that it also does not co-divide, which follows from
    the observations above.
\end{prf}

In the arguments above nonsplitting can be replaced with
Lascar-nonsplitting, and a model $M$ with a set $A$. We can conclude
the following quite desirable characterization of strict nonforking:

\begin{cor}
\begin{enumerate}
Let $N$ be saturated enough over $A$. Then
\item
    A type $p \in \tS(N)$ is
    strictly free over $A$ if and only if for every $a\models p$
    \begin{itemize}
    \item
        $a \ind_A N$
    \item
        $\tp(N/Aa)$ does not divide over $A$.
    \end{itemize}
\item
    If $p \in \tS(N)$ is a heir of $p\rest A$ and
    does not fork over $A$, it is strictly free over $A$. In particular
    this is the case if $p$ is both a heir and a co-heir of $p\rest
    A$.
\end{enumerate}
\end{cor}
\begin{prf}
    (i) Follows directly from the previous Corollary.
    For (ii), note that
    if $p$ is a heir of $p \rest A$, then for every $a \models p$ the type
    $\tp(N/Aa)$ is finitely satisfiable in $A$ (hence does not
    divide  over $A$).
\end{prf}

\begin{dfn}
\begin{enumerate}
\item
    Let $O$ a linear order, $A$ a set. We call a sequence $I =
    \inseq{a}{i}{O}$ a \emph{strict Morley sequence over $B$ based on $A$} if it is
    an indiscernible sequence over $B$ and
    $\tp(a_i/Ba_{<i})$ is strictly free over $A$ for all $i \in O$.
\item
    In the previous definition, we omit ``based on $A$'' if $A=B$.
\item
    Let $p \in \tS(B)$ be a type. We call a sequence $I$
    a \emph{strict Morley sequence in} $p$ if it is a strict Morley
    sequence over $B$ of realizations of $p$.
\end{enumerate}
\end{dfn}

\begin{dfn}
\begin{enumerate}
\item
    We call a type $p \in \tS(A)$ \emph{strictly extendible} if there
    exists a global type extending $p$ which is strictly free over $A$.
\item
    We call a set $A$ a \emph{strict extension base} if every type
    over $A$ is strictly extendible.
\end{enumerate}
\end{dfn}

In \cite{OnUs1} Alf Onshuus and the author show a weak version of
Kim's Lemma, that is: if $\ph(x,a)$ divides over a set $A$ and
witnessed by a sequence $J$ and $\tp(J/A)$ is extendible, then there
is a Morley sequence $I$ in $\tp(a/A)$ which witnessed dividing,
that is, the set $\ph(x,I)$ is inconsistent. Having one such Morley
sequence is often insufficient for applications, though, and we were
wondering whether for some stronger notion of a Morley sequence
``there is'' above can be replaced with ``for all''. It turns out
that strict nonforking provides us with exactly what we need.

\begin{prp}\label{prp:kim}
    Suppose that $\ph(x,a)$ divides over a set $A$ witnessed by an $A$-indiscernible sequence $J$ such
    that $\tp(J/A)$ is extendible. Assume furthermore that $\tp(a/A)$ is
    strictly extendible. Then every strict Morley sequence $I$
    in $\tp(a/A)$ witnesses dividing, that is, the set
    $$\ph(x,I) = \set{\ph(x,a')\colon a' \in I}$$ is inconsistent.
\end{prp}
\begin{rmk}
    Note that since $\tp(a/A)$ is strictly extendible, there exist
    strict Morley sequences in this type.
\end{rmk}
\begin{prf}
    Let $p$ be a global type
    extending $\tp(a/A)$ which is strictly free over $A$. Let $M$ be
    an $(|A|+|T|)^+$-saturated model containing $A$. Let $a' \models
    p \rest M$. Without loss of generality (by applying an
    automorphism), $a' = a$ and $J$ starts with $a$. Denote $J =
    \lseq{a}{\al}{\om}$.

    Since $p$ is strictly free over $A$,  $\tp(M/Aa)$ does not divide over $A$.
    Since $J$ is indiscernible over $A$, it follows that there
    exists $J' \equiv_{Aa} J$ which is indiscernible over $N$.


    Let $Q = \tp(J/M) = Q(x_0,x_1,x_2,\ldots)$. Construct a sequence $\lseq{J}{i}{\om}$ in
    $M$ by $J_i \models Q \rest A\cup J_{<i}$.  Denote $J_i =
    \seq{a_{i,\al} \colon\al<\om}$.

    Note that
    \begin{itemize}
    \item[(*)]
        There is a unique type of ``an infinite sequence in $p$ over
        $A$''. In other words, let (for $\ell = 1,2$) $I_\ell =
        \seq{b_{\ell,\be\colon \be<\om}}$ be such that $b_{\ell,\be}
        \models p\rest Ab_{\ell,<\be}$. Then $\tp(I_1/A) =
        \tp(I_2/A)$. Let us call this type $P$.
    \item[(**)]
        For any $I_1,I_2 \models P$, the set $\ph(x,I_1)$ is
        consistent if and only if $\ph(x,I_1)$ is.
    \item[(***)]
        Let $\eta\colon\omega\to\omega$. Since $J$ is an indiscernible sequence over $M$ of realizations
       of $p\rest M$, one easily sees that the sequence
       $I_\eta = \seq{a_{i,\eta(i)}\colon i<\om}$ realizes $P$.
    \end{itemize}

    So we can conclude

    \begin{itemize}
    \item[$\blacklozenge$]
        For any $I \models P$, the set $\ph(x,I)$ is inconsistent.
    \end{itemize}

    If not, by (**), (***) above we see that
   \begin{itemize}
    \item[$\circ$]
        $\ph(x,J_i)$ is $k$-inconsistent for all $i$
    \item[$\circ\circ$]
        For any $\eta\colon\omega\to\omega$, the set $\ph(x,I_\eta)$
        is consistent.
    \end{itemize}

    This is an easy contradiction to dependence. The most
    straightforward argument would be that this gives the tree
    property of the second kind, $TP_2$, which implies the
    independence property. Since we haven't defined $TP_2$ here, let
    us give some details. By compactness and an Erd\"{o}s-Rado argument, we may assume that
    the sequences $J_i$ are mutually $A$-indiscernible, and $\otp(J_i) = \setQ$.
    We now claim that for every $\eta\colon\omega\to\omega$ the set
    $$\{\ph(x,b_{i,\eta(i)})\}\cup \{\neg \ph(x,b_{i,\al})\}_{\al\neq \eta(i)}$$
    is consistent. Indeed, if $d \models \ph(x,I_\eta)$, then since
    $\ph(x,J_i)$ is $k$-inconsistent, we have that
    $\neg\ph(x,b_{i,\al})$ for almost all $\al$. By mutual
    indiscernibility (and choosing an infinite subsequence of
    $J_i$), the consistence of the set above follows. This
    contradicts dependence.

    \medskip

    So we have shown $\blacklozenge$ above, which states that if $p$ is a strictly free global extension
    of $\tp(a/A)$, and $P$ is constructed as in (*) above, then any
    $I \models P$
    exemplifies dividing of $\ph(x,a)$.

    Now let $I$ be \emph{any} strict Morley sequence in $\tp(a/A)$, and let $p$
    be a global extension of $\Av(I,A\cup I)$ which is strictly free
    over $A$. Constructing $P$ as in (*) above, clearly $I \models
    P$. So $I$ exemplifies dividing of $\ph(x,a)$, as required.


\end{prf}

We have assumed in the Proposition above that a certain type is
strictly extendible. It would be nice to know that there are
``enough'' strictly extendible types. We will show that any type
over a model is such. In fact, we show more: we prove that any type
$p$ over a model $M$ has a global extension $q$ which is a
nonforking heir, that is, $q$ is a heir of $p$ which does not fork
over $M$. The proof relies on the following result due to Chernikov
and Kaplan \cite{ChKa}:

\begin{fct}\label{fct:broom}
    Let $M$ be a model, $p \in \tS(M)$. Assume that
    $$p \vdash \bigvee_{i<k} \ph_i(x,b_i)\lor\bigvee_{j<n}\ps_j(x,c_j)$$
    where $\ph_i(x,y_i), \psi_i(x,z_j)$ are over $M$, $\ph_i(x,b_i)$ does not divide over $M$ for all $i$, and
    $\psi_j(x,c_j)$ divides over $M$ for all $j$. Then there are
    $m<\om$ and automorphisms $\sigma_0, \ldots, \sigma_{m-1}$ over $M$ such that
    $$p \vdash \bigvee_{i<k} \bigvee_{\ell<m}\ph_i(x,\sigma_\ell(b_i))$$
\end{fct}

\begin{cor}\label{cor:heir}
    Let $M$ be a model, $p \in \tS(M)$. Then there exists a global
    heir of $p$ which does not fork over $M$.
\end{cor}
\begin{prf}
    Assume not, then
    $$p \vdash \bigvee_i \ph_i(x,b_i)\lor\bigvee_j\ps_j(x,c_j)$$
    where $\ph_i(x,y_i), \psi_i(x,z_j)$ are over $M$, $\neg\ph_i(x,m) \in p$ for every $m \in M$ for all $i$, and
    $\psi_j(x,c_j)$ divides over $M$ for all $j$. By Fact \ref{fct:broom},
    there is are automorphisms $\sigma_0, \ldots, \sigma_{m-1}$ over
    $M$ such that
    $$p \vdash \bigvee_i \bigvee_{\ell<m}\ph_i(x,\sigma_\ell(b_i))$$

    Continuing as in the well-known proof of existence of heirs, we
    get for some $\theta(x) \in p$
    $$ \models\exists\y\forall x [\theta(x) \rightarrow \bigvee_i \bigvee_{\ell<m}\ph_i(x,\sigma_\ell(b_i))]$$
    hence
    $$ M \models\exists\y\forall x [\theta(x) \rightarrow \bigvee_i \bigvee_{\ell<m}\ph_i(x,\sigma_\ell(b_i))]$$
    which is clearly impossible.
\end{prf}

\begin{rmk}
    We have just shown that any model is a strong
    extension base.
\end{rmk}


An easy conclusion is a particular case of Chernikov and Kaplan's
theorem \cite{ChKa}: in a dependent theory, dividing and forking
coincide over a model.

\begin{cor}(Chernikov,Kaplan, \cite{ChKa})\label{forkingdividing}
    Let $\ph(x,a)$ be a formula, $M$ a model. Then $\ph(x,a)$
    divides over $M$ if and only if $\ph(x,a)$ forks over $M$.
\end{cor}
\begin{prf}
    Suppose $\ph(x,a)$ forks over $M$; so
    $$\ph(x,a) \vdash \bigvee_{i<k}\ps_i(x,a_i)$$
    where each $\ps_i(x,a_i)$ divides over $M$. Assume without loss
    of generality that $\ph(x,a) =
    \bigvee_{i<k}\ps_i(x,a_i)$.

    Let $I = \seq{a_{<k}^\al\colon\al<\om}$ be a strict Morley
    sequence in $\tp(a_{<k}/M)$ (exists since $M$ is a strict
    extension base). It is enough to show that $\ph(x,I)$ is
    inconsistent. Suppose not, and let $b\models \ph(x,I)$. Then
    for some $i<k$, we have $b \models \ps_i(x,a^\al_i)$ for
    infinitely many $\al<\om$. But the sequence $I_i = \seq{a^\al_i\colon
    \al<\om}$ is a strict Morley sequence over $M$ starting with
    $a_i$, hence by ``Kim's Lemma'', Proposition \ref{prp:kim}, it
    has to exemplify dividing of $\ps_i(x,a_i)$, hence
    $\ps_i(x,I_i)$ is $k_i$-inconsistent for some $k_i<\om$. This
    gives the desired contradiction.
\end{prf}

%
%


\section{Bounded weight and orthogonality}

In this section we will point out another important property of
strictly nonforking sequences. So we do not study Morley sequences,
but rather nonforking sequences and sets, not necessarily
indiscernible. In order to develop some of their properties, we will
need to understand  collections of indiscernible sequences whose
first elements are sufficiently independent.

We begin with the following definition:

\begin{dfn}
    Let $A$ be a set and  \lseq{I}{i}{\al} a sequence of sequences. We say that sequences \lseq{I}{i}{\al} are
    \emph{half-mutually $A$-indiscernible} if $I_i$ is indiscernible over
    $AI_{< i}a_{>i}$.
\end{dfn}

\begin{fct}\label{fct:mutind} (Shelah)
    Let $\lseq{a}{i}{\al}$ be a strictly nonforking sequence over $A$, that
    is, $a_i \ind_A^{st} a_{<i}$, and let $I_i$ be an
    $A$-indiscernible sequence starting with $a_i$. Then there exist
    $I'_i \equiv_{Aa_i} I_i$ such that $I'_i$ is indiscernible over
    $AI'_{< i}a_{>i}$ (so \lseq{I'}{i}{\al} are
    half-mutually $A$-indiscernible).
\end{fct}
\begin{prf}
    This is included in \cite{Sh783}, Claim 5.13, but let us still sketch the
    proof. We prove this by induction on $\al$. It is enough to take care of
    the case $\al<\om$.

    Suppose $\lseq{a}{i}{\al+1}$, \lseq{I}{i}{\al+1} are given. By
    the induction hypothesis we may assume that \lseq{I}{i}{\al} are
    half-mutually $A$-indiscernible.

    Since $a_\al \ind^{st}_Aa_{<\al}$, we may assume without loss of
    generality that $a_\al \ind^{st}_AI_{<\al}$. In particular,
    $a_\al \ind_AI_{<\al}$, hence for every $j<\al$ we have $a_\al \ind_{AI_{<j}a_{>j}}I_j$. By preservation of
    indiscernibility, e.g. Observation 8.9 of \cite{Us}, this implies
    that for every $j<\al$ we have $I_{j}$ are indiscernible over $AI_{<j}a_{\le \al, \neq j}$.

    Since $a_\al \ind^{st}_AI_{<\al}$, it also the case that the type
    $\tp(I_{<\al}/Aa_{\al})$ does not divide over $A$. Since $I_\al$ starts with $a_\al$, there is
    $I'_\al$  such that
    \begin{itemize}
    \item
        $I'_{\al} \equiv_{Aa_\al} I_\al$
    \item
        $I'_\al$ is indiscernible over $AI_{<\al}$
    \end{itemize}
    This completes the induction step.

%
%

\end{prf}

The following Corollary is somewhat close to \cite{Sh783}, Claim
5.19, but since we do not think the proof there works as written
(and it is not quite clear whether the claim is true as stated), we
decided to include a precise statement and a proof.

\begin{cor} (Weak Local Character)
\begin{enumerate}
\item
    Let $I = \lseq{a}{i}{|T|^+}$ be a strictly nonforking sequence over $A$ (that is, $a_i \ind^{st}_Aa_{<i}$), $b$
    a finite tuple (or even of cardinality $\le|T|$). Then for
    almost all $i<|T|^+$ (that is, except $|T|$-many) we have that
    \begin{itemize}
    \item
        $\tp(b/Aa_i)$ does not divide over $A$
    \end{itemize}
\item
    If $T$ is strongly dependent, $I$ is an infinite strictly
    nonforking sequence and $b$ is a finite tuple, then for almost
    all (all but finitely many) $a_i \in I$ we have that $\tp(b/Aa_i)$ does not divide over
    $A$.
\end{enumerate}
\end{cor}
\begin{prf}
\begin{enumerate}
\item
    Suppose not. So without loss of generality $tp(b/Aa_i)$
    divides over $A$ for all $i<|T|^+$, and let $I_i$ be an
    $A$-indiscernible sequence starting with $a_i$ witnessing this
    dividing. More precisely, there are formulae $\ph_i(x,y_i)$
    such that
    \begin{itemize}
    \item
        $\ph_i(b,a_i)$ holds
    \item
        The set $\ph(x,I_i) = \set{\ph(x,a')\colon a' \in I_i}$ is
        inconsistent, and even $k_i$-inconsistent for some $k_i$.
    \end{itemize}
    Without loss of generality, $\ph_i(x,y_i) = \ph(x,y)$ and $k_i
    = k$ for all $i<|T|^+$.

    By Fact \ref{fct:mutind} there are $I'_i$ such that
    \begin{itemize}
    \item
        $I'_i \equiv_{Aa_i} I_i$
    \item
        $I'_i$ are half-mutually $A$-indiscernible.
    \end{itemize}

    Note that $\ph(x,I_i)$ are $k$-inconsistent for all $i$.
    Let $I_i = \seq{a_{i,j}\colon j<\om}$. We are going to show the following:

    \begin{itemize}
    \item
        For every $\eta\colon |T|^+\to\omega$ the set
        $\set{\ph(x,a_{i,\eta(i)})\colon i<|T|^+}$ is consistent.
    \end{itemize}

    This will contradict dependence as in the proof of ``Kim's
    Lemma'' (Proposition \ref{prp:kim}). In other words, this
    gives the tree property of the second kind, $TP_2$, which
    implies the independence property.

    In fact, we will show that
    \begin{itemize}
    \item
        For every $\eta\colon |T|^+\to\omega$ we have
        $$ \seq{a_{i,\eta(i)}\colon i<|T|^+} \equiv_A
        \lseq{a}{i}{|T|^+}$$
    \end{itemize}
    This will certainly suffice, because the set $\set{\ph(x,a_{i})\colon
    i<|T|^+}$ is consistent (witnessed by $b$).

    We prove by induction on $\al<|T|^+$ that
    $$ \seq{a_{i,\eta(i)}\colon i<\al}^\frown\seq{a_i\colon \al\le i<|T|^+} \equiv_A
        \lseq{a}{i}{|T|^+}$$

    The case $\al=0$ is trivial, and
    for limit stages use compactness. So let us take care of a
    successor stage. Assume $ \seq{a_{i,\eta(i)}\colon i<\al}^\frown\seq{a_i\colon \al\le i<|T|^+} \equiv_A
        \lseq{a}{i}{|T|^+}$. Since $I_\al$ is indiscernible over
        $AI_{<\al}a_{>\al}$, we have
        $$ \seq{a_{i,\eta(i)}\colon i<\al}^\frown\seq{a_i\colon \al\le i<|T|^+} \equiv_A
        \seq{a_{i,\eta(i)}\colon i<\al}^\frown\seq{a_{\al,\eta(\al)}}\frown\seq{a_i\colon \al <
        i<|T|^+}$$
        which finishes the proof.
\item
    The proof is very similar. That is, using the same arguments we
    arrive at the following situation: for $i<\om$, there are formulae
    $\ph_i(x,a_i)$ which divide over $A$ as exemplified by the
    sequences $I_i$, whereas for every $\eta\colon\omega\to\omega$
    the set $\set{\ph(x,a_{i,\eta(i)})\colon i<|T|^+}$ is
    consistent. As in the proof of Proposition \ref{prp:kim}, we
    may assume by compactness that $I_i$ are of order type
    $\setQ$, so $I_i = \seq{a_{i,q}\colon q\in \setQ}$, with $a_i =
    a_{i,0}$. It is now easy to see (by taking infinite subsequences of $I_i$, since $\ph(x,I_i)$
    is $k_i$-inconsistent for some $k_i$) that the following set is
    consistent for every $\eta\colon\omega\to\omega$:
    $$\{\ph_i(x,b_{i,\eta(i)})\colon i\in\setQ\}\cup \{\neg \ph_i(x,b_{i,\al})\colon i\in\setQ, \al\neq \eta(i)\}$$

    And this is precisely the definition of lack of strong
    dependence.
\end{enumerate}
\end{prf}

Note that one can regard this ``weak local character'' as a kind
of ``bounded pre-weight'', or ``rudimentarily finite pre-weight''
in the case of strongly dependent theories (developing further
some concepts introduced by Alf Onshuus and the author in
\cite{OnUs1}). Indeed, we have shown that given an ``independent
enough'' sequence \lseq{a}{i}{|T|^+} and a tuple $b$, it is the
case that $b$ can only divide with a few $a_i$'s (when working
over a model or any set which is a strong extension base, one can
replace dividing with forking by \cite{ChKa}, see Corollary
\ref{forkingdividing}).

There are several natural questions that arise in this context.
For example, given an ordinal $\al$, one can define $p = tp(b/A)$
to have \emph{forking pre-weight} at least $\al$ if there are
$\set{a_i\colon i<\al}$ forking independent over $A$ (that is,
$a_i \ind_A a_{\neq i}$) and $\tp(b/Aa_i)$ divides over $A$ for
all $i$. On the other hand, we can define $p$ to have \emph{strict
forking pre-weight} at least $\al$ if there is $\seq{a_i\colon
i<\al}$ a strictly forking independent sequence over $A$ (that is,
$a_i \ind^{st}_A a_{< i}$) and $\tp(b/Aa_i)$ divides over $A$ for
all $i$. Let us say that $p$ has \emph{rudimentarily finite}
pre-weight (forking, strict forking, etc) if it is not the case
that the pre-weight of $p$ is at least $\omega$.

\begin{cor}
\begin{enumerate}
\item
    If $T$ is strongly dependent, then every type has rudimentarily finite
    strict forking pre-weight. In fact, if a type $p$ is strongly
    dependent (as defined in \cite{OnUs1}, Definition 2.6), then $p$ has rudimentarily
    finite strict forking pre-weight.
\item
    If every type in $T$ has rudimentarily finite forking
    pre-weight, then $T$ is strongly dependent. In fact, if a type
    $p$ has rudimentarily finite forking
    pre-weight, then $p$ is strongly dependent.
\end{enumerate}
\end{cor}
\begin{prf}
\begin{enumerate}
    \item By the previous Corollary.
    \item This is a trivial consequence of \cite{OnUs1}, Theorem
    2.12(ii).
\end{enumerate}
\end{prf}

We see that: $T$ has rudimentarily finite forking weight
$\implies$ $T$ is strongly dependent $\Longleftrightarrow$ $T$ has
rudimentarily finite weight in the sense of \cite{OnUs1},
Definition 2.3 $\implies$ $T$ has rudimentarily finite strict
forking weight. So it is natural to wonder

\begin{qst}
    Are any of the above implications reversible?
\end{qst}

Reading the proofs of the results in this section carefully, one
sees that the main issue has to do with appropriate notions of
\emph{weak orthogonality}. Let us give several possible
definitions and point out some connections between them.

\begin{dfn}
\begin{enumerate}
    \item (Shelah, e.g. \cite{Sh783}, Definition 5.32) We call two
    types $p,q \in \tS(A)$ \emph{weakly orthogonal} or if $p(x)\cup
    q(y)$ is a complete type over $A$. We write $p \perp_w q$. If
    $a,b$ realize $p,q$ respectively, then we write $a/A\perp_w
    b/A$ or $a\perp_wb$ when $A$ is fixed and clear from the
    context.
    \item We call $\tp(a/A), \tp(b/A)$ \emph{weakly
    orthogonal$^{1}$} if whenever $I,J$ are $A$-indiscernible
    sequences starting with $a,b$ respectively, there are $I',J'$
    mutually $A$-indiscernible such that $I \equiv_{Aa} I'$ and $J
    \equiv_{Ab} J'$. We write $a/A\perp^1_w
    b/A$ or $a\perp^1_wb$ when $A$ is fixed and clear from the
    context.
    \item We call $\tp(a/A), \tp(b/A)$ \emph{weakly
    orthogonal$^{\frac 12}$} if whenever $I,J$ are $A$-indiscernible
    sequences starting with $a,b$ respectively, there are $I',J'$
    half-mutually $A$-indiscernible such that $I \equiv_{Aa} I'$ and $J
    \equiv_{Ab} J'$. We write $a/A\perp^{\frac 12}_w
    b/A$ or $a\perp^{\frac 12}_wb$.
    \item We call $\tp(a/A), \tp(b/A)$ \emph{weakly
    orthogonal$^{st}$} if whenever $a\models p$ and $b
    \models q$, we have $a \ind^{st}_A b$ and $b \ind^{st}_A a$. We write $a/A\perp^{st}_w
    b/A$ or $a\perp^{st}_wb$.
    \item We call $\tp(a/A), \tp(b/A)$ \emph{weakly
    orthogonal$^{fk}$} if whenever $a\models p$ and $b
    \models q$, we have $a \ind_A b$ and $b \ind_A a$. We write $a/A\perp^{fk}_w
    b/A$ or $a\perp^{fk}_wb$.
\end{enumerate}
\end{dfn}

\begin{obs}\label{obs:orth} Let $A$ be a set, $p,q \in \tS(A)$.
\begin{enumerate}
  \item If $A$ is an extension base (or just $p,q$ do not fork over $A$) and
  $p\perp_wq$, then $p\perp^{fk}_wq$.
  \item If $A$ is a strict extension base (or just $p,q$ are strictly free over $A$) and
  $p\perp_wq$, then $p\perp^{st}_wq$.
  \item If $p\perp^{st}_wq$ then $p\perp^{fk}_wq$.
  \item If $p\perp^{st}_wq$ then $p\perp^{\frac 12}_wq$.
  \item If $p\perp^{1}_wq$ then $p\perp^{\frac 12}_wq$.
\end{enumerate}
\end{obs}
\begin{prf}
    Easy, e.g. (iv) is Fact \ref{fct:mutind}.
\end{prf}

\begin{lem}
    Let $A$ be an extension base(e.g. a model), $p,q \in \tS(A)$. If $p\perp_wq$, then $p\perp^1_wq$.
\end{lem}
\begin{prf}
    Let $I,J$ be $A$-indiscernible starting with $a,b$ respectively
    where $a\models p$ and $b \models q$. Let $I =
    \lseq{a}{i}{\om}$ and $J = \lseq{b}{i}{\om}$.

    By an Erd\"{o}s-Rado argument, there is $I' \equiv_A I$ which is
    indiscernible over $AJ$. Similarly, there is $J' \equiv_A J$
    which is indiscernible over $AI'$, moreover, denoting $J' =
    \lseq{b'}{i}{\om}$, we have that
    \begin{itemize}
    \item[$(\blacklozenge)$]
    for every $k<\om$ there are
    $j_1, \ldots j_k$ such that $\lseq{b'}{i}{k} \equiv_{AI'}
    \seq{b_{j_i}\colon i<k}$.
    \end{itemize}

    Denote $I' = \lseq{a'}{i}{\om}$. We claim that it is still the case that
    $I'$ is indiscernible over $AJ'$. Suppose not; so there
    is a formula $\ph(\x,\b')$ over $AJ'$ such that
    $\ph(a'_{i_1},\ldots,a'_{i_k},\b')\land\neg\ph(a'_{j_1},\ldots,a'_{j_k},\b')$
    for some $i_1<\ldots<i_1<\om$ and $j_1<\ldots<j_1<\om$. But by $(\blacklozenge)$
    above,
    there is a tuple $\b$ of elements of $J$ satisfying the same
    formula, that is,
    $\ph(a'_{i_1},\ldots,a'_{i_k},\b)\land\neg\ph(a'_{j_1},\ldots,a'_{j_k},\b)$,
    which implies that $I'$ is not indiscernible over $AJ$, a
    contradiction.

    Finally, let us note that $I',J'$ have the same type over $A$ as
    $I,J$ respectively. Hence in particular $a'_0 \models p$ and
    $b'_0 \models q$. By the assumption $p\perp_wq$, we have
    $a'_0b'_0 \equiv_A ab$. So without loss of generality $a'_0b'_0
    = ab$, and $I',J'$ are as required in the definition of
    $p\perp^1_wq$.

\end{prf}

\begin{exm}
    Let $(\setQ,<,P)$ be the theory of $(\setQ,<)$ with a dense
    co-dense predicate $P$. Let $p,q \in \tS(\setQ)$ be the types
    over the prime model $\setQ$ such that if $a,b$ realize $p,q$
    respectively, then $a,b > \setQ$, $P(a), \neg P(b)$.

    Clearly $p \not\perp_w q$ since $p,q$ do not determine whether
    $a<b$ or $b<a$. On the other hand, it is easy to see that $p
    \perp^1_wq$ and $p \perp^{st}_wq$.

    Hence the implications in the Lemma above, as well as Observation \ref{obs:orth}(i),(ii) are not reversible.
\end{exm}

So we obtain:

\begin{cor}
    Let $A$ be a strict extension base (e.g. a model), $p,q \in \tS(A)$.
    Then
    $$p \perp_w q \implies p\perp^{st} q \implies p \perp^{\frac
    12}
    q$$
    $$p \perp_w q \implies p\perp^{1} q \implies p \perp^{\frac
    12}
    q$$
    $$p \perp_w q \implies p\perp^{st} q \implies p \perp^{fk}$$
\end{cor}

\bigskip

Again, there are many natural questions.

\begin{qst}\label{qst:orth}
\begin{enumerate}
  \item Which of the implications above are reversible? We know that
  the first one in each row is \emph{not}.
  \item What is the relation between $\perp^1_w$ and $\perp^{st}_w$?
  \item What is the relation between $\perp^{\frac 12}_w$ and $\perp^{fk}_w$?
  \item More specifically, is it the case that whenever $a \ind_Ab$
  and $b\ind_Aa$ (and $A$ is sufficiently nice), then $a
  \ind^{st}b$? The reverse is clearly true.
\end{enumerate}
\end{qst}

The following is a small step in the direction of (possibly)
answering Question \ref{qst:orth} (iii) or (iv):

\begin{lem}
    Let $A$ be a set, $I_i$ (for $i<k$) be a sequence starting with the element
    $a_i$ such that $I_i$ is indiscernible over $A$. Assume furthermore that
    $\set{a_i\colon i<k}$ is forking independent over $A$. Then
    without loss of generality $I_i$ is indiscernible over $Aa_{\neq i}$.
    In other words, there are $I'_i$ such that
    \begin{itemize}
    \item
        $I'_i \equiv_{Aa_i} I_i$
    \item
        $I'_i$ is indiscernible over $Aa_{\neq i}$
    \end{itemize}
\end{lem}
\begin{prf}
    We prove this by induction on $k$, the case $k=1$ being trivial.

    So let $k>1$ and assume that $I_i$ is indiscernible over
    $Aa_{<k, \neq i}$ for all $i<k$.

    Recall that $a_{<k}\ind_Aa_k$. Let $a'_{<k} \equiv_{Aa_k}
    a_{<k}$ be such that $a'_{<k} \ind_A I_k$. Let $\sigma \in
    \Aut(\fC/Aa_k)$ take $a'_{<k}$ to $a_{<k}$, and denote $I'_k = \sigma(I_k)$. Clearly
    \begin{itemize}
    \item
        $I'_k \equiv_{Aa_k} I_k$
    \item
        $a_{<k} \ind_A I'_k$
    \end{itemize}
    Hence by preservation of indiscernibility (e.g. Observation 8.9 in \cite{Us}), $I'_k$ is indiscernible over $Aa_{<k}$.

    Now since $a_k \ind_A a_{<k}$, we can find (as before) $I'_i$
    for $i<k$ satisfying

   \begin{itemize}
    \item
        $I'_{<i} \equiv_{Aa_{<k}} I_i$
    \item
        $a_{k} \ind_A I'_{<k}$
    \end{itemize}
\emph{}
    In particular, we obtain for every $i<k$:

    \begin{itemize}
    \item
        $a_k \ind_A a_{<k}I_i$, hence $a_k \ind_{Aa_{<k,\neq
        i}}I'_i$. By the induction hypothesis and preservation of indiscernibility this implies
        that $I'_i$ is indiscernible over $Aa_{\neq i}$.
    \end{itemize}

    Recall that $I'_k$ is indiscernible over $Aa_{<k}$. So we are
    done.

\end{prf}

One of the earlier versions of \cite{Sh783} contained the statement
that the conclusion of the Lemma above can be obtained when starting
from a weaker assumption: the sequence $\lseq{a}{i}{k}$ is
nonforking. The following example shows that this is not always
possible:

\begin{exm}
    Consider the theory of $(\setQ,<)$, and let $b = \seq{0,2}, a =
    1$. Then $b \ind a$, but if $I =
    \left<\seq{0,2},\seq{3,5},\seq{6,8}, \ldots\right>$ and $J =
    \seq{1,4,7,\ldots}$, then clearly there are no $I',J'$ of the
    same type as $I,J$ respectively, starting with $a,b$ such that
    $J'$ is indiscernible over $a$.

    We constructed the example with $A = \emptyset$, but it is as
    easy to modify it such that $A$ is any set, in particular a model.
\end{exm}

Let us conclude this section with the following Lemma, which is a
slight generalization of Proposition 4.13 in \cite{OnUs1}, and
ideas behind the proof are very similar. We include it because it
might also become useful for questions related to issues discussed
above (although right now we do not see any concrete
applications).

\begin{lem}
    Let $A$ be a set, $I_i$ (for $i<k$) be a sequence starting with the element
    $a_i$ such that:
    \begin{itemize}
    \item[($\blacklozenge$)]
        $I_i \ind_A a_{<i}$
    \item[($\blacklozenge\blacklozenge$)]
        $I_i$ is indiscernible over $Aa_{< i}$
    \end{itemize}
    Then
    without loss of generality $I_i$ is indiscernible over $AI_{\neq
    i}$. Moreover, there are $I'_i$ such that
    \begin{itemize}
    \item
        $I'_i \equiv_{Aa_i} I_i$
    \item
        $I'_i$ is indiscernible over $AI'_{\neq i}$
    \item
        $I'_i \ind_A I'_{<i}$
    \end{itemize}
\end{lem}
\begin{prf}
    We prove the lemma by induction on $k$, the case $k=1$ being
    trivial.

    Let $a_{\le k}$, $I_{\le k}$ be as in the assumptions of the lemma. By
    the induction hypothesis we can assume that for $i<k$ we have
    sequences $I''_i$ satisfying the conclusion, that is

    \begin{enumerate}
    \item
        $I''_i \equiv_{Aa_i} I_i$
    \item
        $I''_i$ is indiscernible over $AI'_{\neq i}$
    \item
        $I''_i \ind_A I''_{<i}$
    \end{enumerate}

    So in particular each $I''_i$ (for $i<k$) starts with the element
    $a_i$.

    Recall that by the assumptions on $I_k$ we also have
    \begin{itemize}
    \item[(iv)]
        $I_k$ is indiscernible over $Aa_{<k}$
    \item[(v)]
        $I_k \ind_A a_{<k}$
    \end{itemize}

    By (v) above and the existence of nonforking extensions, there
    are $I^*_{<k}$ such that

    \begin{itemize}
    \item
        $I^*_{<k} \equiv_{Aa_{<k}} I''_{<k}$
    \item
        $I_k \ind_A I^*_{<k}$
    \end{itemize}

    So without loss of generality we may assume in addition that
    \begin{itemize}
    \item[(vi)]
        $I_k \ind_A I''_{<k}$
    \end{itemize}

    Since we can make $I_k$ as long as we wish, applying Erd\"{o}s-Rado, there exists $I''_k$ such that

    \begin{itemize}
    \item[(*)]
        $I''_k$ is indiscernible over $AI''_{<k}$
    \item[(**)]
        Every $n$-type of $I''_k$ over $AI''_{<k}$ ``appears'' in
        $I_k$
    \end{itemize}

    Note:

    \begin{itemize}
    \item[($\diamondsuit$)]
        $I''_k \equiv_{Aa_{<k}} I_k$ [since $I_k$ satisfies
        (iv) above and $I''_k$ satisfies (*),(**)]
    \item[($\diamondsuit\diamondsuit$)]
        $I''_k \ind_A I''_{<k}$ [by (vi) and (**) above]
    \end{itemize}

    Let $\sigma \in \Aut(\fC/Aa_{<k})$ be such that
    $\sigma(I''_k) = I_k$. Define $I'_i = \sigma(I''_i)$ for
    $i<k$ and $I'_k = I_k$. We claim that $I'_{\le k}$ satisfy the
    conclusion of the lemma, which completes the induction step.
    Indeed,

    \begin{itemize}
    \item
        $I'_i \equiv_{Aa_i} I_i$ for all $i \le k$:
        \begin{itemize}
        \item $i<k$. By the
            induction hypothesis + $\sigma$ being over $a_{<k}$.
        \item $i = k$. Clear since $I'_k = I_k$.
        \end{itemize}

    \item
        For every $i\le k$ we have that $I_i$ is indiscernible over $AI_{\neq
        i}$:
        \begin{itemize}
        \item
            $i<k$.  By
        $(\diamondsuit\diamondsuit)$ above we have $I'_k \ind_A I'_{<k}$, hence
        $I'_k \ind_{B} I'_i$ where $B = AI'_{<k,\neq i}$. By the induction hypothesis,
        $I'_i$ is indiscernible over $B$, hence by preservation of indiscernibility, $I'_i$ is indiscernible over
        $AI'_{\neq i}$, as required.
        \item
            $i=k$. Follows immediately from the choice of $I''_k$.
        \end{itemize}

    \item
        $I'_i \ind_A I'_{<i}$ for all $i$:
        \begin{itemize}
        \item $i<k$. By the
            induction hypothesis.
        \item $i = k$. By ($\diamondsuit\diamondsuit$) above.

        \end{itemize}

    \end{itemize}

\end{prf}

\begin{cor}\label{cor:Morleyindisc}
    Let $A$ be a set, $I_i$ (for $i<k$) be a sequence starting with the element
    $a_i$ such that:
    \begin{itemize}
    \item[($\blacklozenge$)]
        $I_i$ is a Morley sequence over $Aa_{<i}$ based on $A$
    \end{itemize}
    Then
  there are $I'_i$ such that
    \begin{itemize}
    \item
        $I'_i \equiv_{Aa_i} I_i$
    \item
        $I'_i$ is indiscernible over $AI'_{\neq i}$
    \item
        $I'_i \ind_A I'_{<i}$
    \end{itemize}
\end{cor}
\begin{prf}
    Follows immediately from the previous lemma (and transitivity of nonforking on
    the left).
\end{prf}

Some of the questions in this section will be addressed and
partially answered in a subsequent work of Itay Kaplan and the
author \cite{KaUs}.

\bibliography{common}
\bibliographystyle{alpha}

\end{document}